\documentclass[11pt]{amsart}
\usepackage[all]{xy}

\usepackage{amsmath}
\usepackage[latin1]{inputenc}
\usepackage{amsfonts}
\usepackage{amssymb}
\usepackage{multicol}
\usepackage{verbatim}
\usepackage{amsthm}
\usepackage{amscd}
\usepackage{pb-diagram}
\usepackage[mathscr]{eucal}
\usepackage[latin1]{inputenc}

\usepackage{amsmath}
\usepackage{epsfig}
\usepackage{amsfonts}

\usepackage{amssymb}

\usepackage{amsthm}

\usepackage{graphics}

\newtheorem{thm}{Theorem}[section]

\newtheorem{defi}[thm]{Definition}

\newtheorem{cor}[thm]{Corollary}

\newtheorem{prop}[thm]{Proposition}

\newtheorem{conj}[thm]{Conjecture}

\newtheorem{rem}[thm]{Remark}

\usepackage[all]{xy}
\title{Advances in R-matrices and their applications [after Maulik-Okounkov, Kang-Kashiwara-Kim-Oh...] }
\author{David HERNANDEZ}\thanks{Supported
in part by the European Research Council under the European Union's Framework Programme
H2020 with ERC Grant Agreement number 647353 QAffine. }
\address{Sorbonne Paris Cit\'e\\
Universit\'e Paris-Diderot Paris 7\\
Institut de Math\'ematiques de Jussieu--Paris Rive Gauche\\
CNRS UMR 7586\\
B\^atiment Sophie Germain\\
F--75205 Paris Cedex 13}
\email{david.hernandez@imj-prg.fr}

\begin{document}

\maketitle

\pagestyle{myheadings}

\markboth{DAVID HERNANDEZ}{ADVANCES IN R-MATRICES AND THEIR APPLICATIONS}

This is an English translation of the Bourbaki Seminar no. 1129 (March 2017).

\bigskip

\noindent{\bf INTRODUCTION}

\bigskip
\bigskip

R-matrices are the solutions of the Yang-Baxter equation. At the origin of the quantum group theory, they may
be interpreted as intertwining operators. Recent advances have been made independently in 
different directions. Maulik-Okounkov have given a geometric approach to R-matrices with new tools
in symplectic geometry, the stable envelopes. Kang-Kashiwara-Kim-Oh proved a conjecture on the 
categorification of cluster algebras using R-matrices in a crucial way. Eventually, a better
understanding of the action of transfer-matrices obtained from R-matrices led to the proof of
several conjectures about the corresponding quantum integrable systems.

\bigskip
\bigskip
\bigskip

The Yang-Baxter equation
$$\mathcal{R}_{12}(z)\mathcal{R}_{13}(zw)\mathcal{R}_{23}(w) = \mathcal{R}_{23}(w)
\mathcal{R}_{13}(zw) \mathcal{R}_{12}(z)$$
is written in terms of a  Laurent power series
$$\mathcal{R}(z)\in (\mathcal{A}\otimes \mathcal{A})((z))$$
with coefficients in the tensor square of an algebra $\mathcal{A}$ that we suppose to be complex. We use the standard notation
$$\mathcal{R}_{12}(z) = \mathcal{R}(z)\otimes 1\text{ , }\mathcal{R}_{23}(z) = 1 \otimes \mathcal{R}(z)\text{ , }
\mathcal{R}_{13}(z) = (P\otimes \text{Id})(\mathcal{R}_{23}(z)),$$
where $P$ is the twist satisfying 
$$P(x\otimes y) = y\otimes x.$$
So the terms in the Yang-Baxter equation belong to $(\mathcal{A}^{\otimes 3})((z,w))$. A solution is called an $R$-matrix, or affine $R$-matrix as it depends here on the parameter $z$, called a spectral (or affine) parameter.

This equation originated in the theory of (quantum) integrable systems. It was also intensively studied in various other fields, in particular in representation theory and in topology (we refer to the seminars \cite{Ro, STS}). 
The theory of quantum groups has been introduced to give answers to the problem of the construction of such $R$-matrices.

\section{Algebraic construction}\label{algcons}

\subsection{Fundamental example}\label{exfund}

Let us start with the following fundamental example for the algebra $\mathcal{A} = \text{End}(V)$ where $V$ is 
a $2$-dimensional complex vector space with a basis $\mathcal{B}$. For $q\in\mathbb{C}^*$, we have the $R$-matrix written in the basis of $V\otimes V$ naturally associated to $\mathcal{B}$ :
\begin{equation}\label{fondsl2}\begin{pmatrix}1&0&0&0\\ 0&\frac{q^{-1}(z-1)}{z-q^{-2}}&\frac{1-q^{-2}}{z - q^{-2}}&0\\0&\frac{z(1 - q^{-2})}{z-q^{-2}}&\frac{q^{-1}(z - 1)}{z - q^{-2}}&0\\0&0&0&1\end{pmatrix}\in (\text{End}(V^{\otimes 2}))((z)) \simeq (\mathcal{A}^{\otimes 2})((z))  .\end{equation}
By considering the submatrix (after having removed the first and last lines and columns), setting $z = e^u$, $q = e^{h/2}$ and for $u$ and $h$ close to $0$, one gets the famous ``Yang $R$-matrix":
\begin{equation}\label{yang}\frac{1}{u+h}\begin{pmatrix} u & h\\ h & u\end{pmatrix}.\end{equation}
Historically, the $R$-matrix (\ref{fondsl2}) appeared in statistical physics in the framework of the study of the 6-vertex model introduced by Pauling (1935) which gives, in particular, a description of ice crystals. The study of this model is closely related to the study of another model, in quantum statistical physics, called the $XXZ$-model of $1/2$ spin or the quantum Heisenberg model (1928). It is a model for quantum magnetic spin chains with two corresponding classical states, up or down.

These two models, the $6$-vertex model and the $XXZ$-model, are some of the most studied models in statistical and quantum
physics. The underlying mathematical structures are closely related and involve the quantum affine algebra 
$\mathcal{U}_q(\hat{\mathfrak{sl}_2})$ associated to the Lie algebra $\mathfrak{g} = \mathfrak{sl}_2$.
This algebra has a family of $2$-dimensional simple representations $V$ called ``fundamental representations".
The $R$-matrix (\ref{fondsl2}) can be obtained from such representations. But the theory of quantum groups produces many more examples, by changing the Lie algebra $\mathfrak{g}$ or the representation $V$, which correspond to many more quantum systems.

 \subsection{Universal $R$-matrices}

The example discussed above can be generalized by using the following procedure to construct large families of $R$-matrices. In the following we suppose that $q \in\mathbb{C}^*$ is not a root of unity. To a complex simple finite-dimensional Lie algebra $\mathfrak{g}$ are associated :

\medskip

- on one hand the affine Kac-Moody algebra $\hat{\mathfrak{g}} = \mathcal{L}(\mathfrak{g}) \oplus \mathbb{C} c$ (see the seminar \cite[Section 3.1]{F}) defined as the universal central extension of the loop algebra  
$$\mathcal{L}(\mathfrak{g}) = \mathfrak{g} \otimes \mathbb{C}[t^{\pm 1}]$$ 

\medskip

- on the other hand the quantum algebra $\mathcal{U}_q(\mathfrak{g})$ (quantum group in the sense of Drinfeld and Jimbo, see the seminar \cite{Ro}) which is a Hopf algebra and a $q$-deformation of the enveloping algebra $\mathcal{U}(\mathfrak{g})$. 

\medskip

If one considers simultaneously these two generalizations of Lie algebras, one obtains the quantum affine algebra 
$\mathcal{U}_q(\hat{\mathfrak{g}})$. It is a Hopf algebra and in particular it has a coproduct. That is, an algebra morphism
$$\Delta : \mathcal{U}_q(\hat{\mathfrak{g}})\rightarrow \mathcal{U}_q(\hat{\mathfrak{g}})\otimes \mathcal{U}_q(\hat{\mathfrak{g}})$$
which gives the structure of a $\mathcal{U}_q(\hat{\mathfrak{g}})$-module on the tensor product of $\mathcal{U}_q(\hat{\mathfrak{g}})$-modules.

Drinfeld \cite{Dri2} established (the proof was detailed subsequently by Beck \cite{bec} and Damiani \cite{da}) that $\mathcal{U}_q(\hat{\mathfrak{g}})$ can be obtained not only as a quantization of $\mathcal{U}(\hat{\mathfrak{g}})$, but also by another process; as an  affinization of the quantum group $\mathcal{U}_q(\mathfrak{g})$. It is the Drinfeld realization of quantum affine algebras. This might be seen in the form of the following commutative diagram :
$$\xymatrix{ &\hat{\mathfrak{g}}\ar@{-->}[dr]^{\text{Quantization}}&   
\\ \mathfrak{g} \ar@{-->}[ur]^{\text{Affinization}}\ar@{-->}[dr]_{\text{Quantization}}& &  \mathcal{U}_q(\hat{\mathfrak{g}})
\\ & \mathcal{U}_q(\mathfrak{g})\ar@{-->}[ur]_{\text{Quantum affinization}}& }$$
Hence we get two distinct presentations of isomorphic algebras (the original Drinfeld-Jimbo presentation and the Drinfeld presentation). It is a quantum analogue of a classical theorem of Kac and Moody for affine Kac-Moody algebras (see \cite[Chapter 7]{ka}).

Using the Drinfeld presentation the algebra $\mathcal{U}_q(\hat{\mathfrak{g}})$ can be endowed with a natural $\mathbb{Z}$-grading. This grading corresponds to automorphisms $\tau_a$ of the algebra $\mathcal{U}_q(\hat{\mathfrak{g}})$ for $a\in\mathbb{C}^*$, as well as to an automorphism $\tau_z$ of the algebra $\mathcal{U}_q(\hat{\mathfrak{g}})(z)$ for an indeterminate $z$ so that for $g\in\mathcal{U}_q(\hat{\mathfrak{g}})$, homogeneous of degree $m\in\mathbb{Z}$, one has :
$$\tau_a(g) = a^m g\text{ and }\tau_z(g) = z^m g.$$

\begin{thm}[Drinfeld] The algebra $\mathcal{U}_q(\hat{\mathfrak{g}})$ has a universal $R$-matrix. That is, a non-trivial solution of the Yang-Baxter equation in the completed tensor product 
$$\mathcal{R}(z)\in [\mathcal{U}_q(\hat{\mathfrak{g}})\hat{\otimes}\mathcal{U}_q(\hat{\mathfrak{g}})][[z]].$$ 
Moreover, for $V$ and $W$ finite-dimensional representations of $\mathcal{U}_q(\hat{\mathfrak{g}})$ the image 
$$(\rho_V\otimes \rho_W)(\mathcal{R}(z)) = \mathcal{R}_{V,W}(z)\in (\text{End}(V)\otimes \text{End}(W))[[z]],$$
for the representation morphisms $\rho_V$ and $\rho_W$, is well-defined and  
$$P\circ \mathcal{R}_{V,W}(z) : (V\otimes W)[[z]]\rightarrow (W\otimes V)[[z]]$$
is a morphism of $\mathcal{U}_q(\hat{\mathfrak{g}})$-modules ($P$ is the twist morphism as above and the action on $V$ is 
twisted by the automorphism $\tau_z$). 

Eventually one has the relations
\begin{equation}\label{proprt}(\text{Id}\otimes \Delta)(\mathcal{R}(z)) = \mathcal{R}_{13}(z)\mathcal{R}_{12}(z)
\text{ and }(\Delta\otimes \text{Id})(\mathcal{R}(z)) = \mathcal{R}_{13}(z)\mathcal{R}_{23}(z).\end{equation}
\end{thm}

Hence one has intertwining operators. When $V = W$, $\mathcal{R}_{V,V}(z)$ is an $R$-matrix.

\begin{rem} The notion of completed tensor product $\hat{\otimes}$ above is defined by using certain filtrations of $\mathcal{U}_q(\hat{\mathfrak{g}})$. We do not give additional details as this notion plays a role only for the universal $R$-matrix and not for the $R$-matrices $\mathcal{R}_{V,W}(z)$ (or for the transfer-matrices) studied in the rest of this seminar.
\end{rem}

\begin{rem} The conventions for the $\mathbb{Z}$-gradation of $\mathcal{U}_q(\hat{\mathfrak{g}})$ discussed above imply that $\mathcal{R}(z)$ has only non-negative powers of $z$.
\end{rem}

\noindent{\sc Proof} (sketch) --- The proof of Drinfeld relies on the realization of $\mathcal{U}_q(\hat{\mathfrak{g}})$ as a quotient of the double of a Hopf subalgebra (analogous to the Borel subalgebra). The double $D(A)$ of a Hopf algebra $A$ is the vector space $A\otimes A^0$ where $A^0$ is the dual Hopf algebra of $A$. Then $D(A)$ is endowed with a Hopf algebra structure which automatically has a universal $R$-matrix (see the seminar \cite[Section 3]{Ro}).\qed

\subsection{Normalized $R$-matrices}\label{rnorm}

We have the following rationality property.

\begin{prop} For $V$ and $W$ finite-dimensional representations of $\mathcal{U}_q(\hat{\mathfrak{g}})$, there exists a unique formal Laurent power series, $f_{V,W}(z)\in\mathbb{C}((z))$, such that the product 
$$f_{V,W}(z)\mathcal{R}_{V,W}(z) \in (\text{End}(V\otimes W))(z)$$
is rational in $z$.

\end{prop}

\noindent{\sc Proof} (sketch) ---  One can use for example the argument in \cite[Section 9.2]{efk} (see also \cite{ifre}).  Let $v$, $w$ be highest weight vectors respectively in $V$ and $W$ 
(relative to an analogue of the Cartan subalgebra in $\mathcal{U}_q(\hat{\mathfrak{g}})$). 
Then there is a unique $f_{V,W}(z)\in\mathbb{C}[[z]]$ so that
\begin{equation}\label{norm}(f_{V,W}(z)\mathcal{R}_{V,W}(z)).(v\otimes w) = v\otimes w.\end{equation}
The representation $(V\otimes W)\otimes \mathbb{C}((z))$ of $\mathcal{U}_q(\mathfrak{g})((z))$ is simple
(this is called the generic irreducibility of the tensor product). Hence the fact that
$f_{V,W}(z)P\circ\mathcal{R}_{V,W}(z)$ is an intertwining operator can be translated into a system of linear equations whose solutions are rational.\qed

\begin{rem}\label{inv} One gets the uniqueness of the Laurent series $f_{V,W}(z)$ by requiring the relation
 (\ref{norm}). This will be assumed in the following. The rational morphism $f_{V,W}(z)\mathcal{R}_{V,W}(z)$ 
has an inverse which is, under this assumption,
$$(f_{V,W}(z)\mathcal{R}_{V,W}(z))^{-1} = f_{W,V}(z^{-1})(P\circ \mathcal{R}_{W,V}(z^{-1})\circ P).$$
\end{rem}

We will be particularly interested in the case of fundamental representations  $V_i(a)$ of $\mathcal{U}_q(\hat{\mathfrak{g}})$. They are parametrized by $i\in\{1,\ldots, n\}$ and $a\in\mathbb{C}^*$ where $n$ is the rank of $\mathfrak{g}$. The reader may refer to \cite[Chapter 12.2]{CP} for general results on finite-dimensional representations of $\mathcal{U}_q(\hat{\mathfrak{g}})$.

Let $R\geq 0$ be the order of $1$ as a pole of $f_{V,W}(z)\mathcal{R}_{V,W}(z)$. Then the limit  
\begin{equation}\label{limr}\text{lim}_{z\rightarrow 1}[(z - 1)^Rf_{V,W}(z)P\circ \mathcal{R}_{V,W}(z)] : V\otimes W\rightarrow W\otimes V\end{equation}
is a non-zero morphism of $\mathcal{U}_q(\hat{\mathfrak{g}})$-module (if $V$ and $W$ are non-zero).

\begin{defi} The limit (\ref{limr}) is called a normalized $R$-matrix and is denoted 
$$\mathcal{R}_{V,W}^{norm} : V\otimes W\rightarrow W\otimes V.$$
\end{defi}

Consider for example $\mathfrak{g} = \mathfrak{sl}_2$ with two fundamental representations $V = V_1(1)$ and $W = V_1(q^2)$ of dimension $2$. In this case they can be constructed from a representation of $\mathcal{U}_q(\mathfrak{sl}_2)$ using evaluation morphisms\footnote{Such evaluation morphisms defined by Jimbo exist only in type $A$. This is a source of important technical difficulties in the study of finite-dimensional representations of $\mathcal{U}_q(\hat{\mathfrak{g}})$.} $\mathcal{U}_q(\hat{\mathfrak{sl}}_2)\rightarrow \mathcal{U}_q(\mathfrak{sl}_2)$. We get the $R$-matrix (\ref{fondsl2}) with $zq^{-2}$ instead of $z$. The pole $1$ is simple and 
$$\text{lim}_{z\rightarrow 1} (z-1) \begin{pmatrix}1&0&0&0\\ 0&\frac{q(zq^{-2}-1)}{z - 1}&\frac{q^2 - 1}{z - 1}&0\\0&\frac{z(1 - q^{-2})}{z- 1}&\frac{q(zq^{-2} - 1)}{z - 1}&0\\0&0&0&1\end{pmatrix} =  \begin{pmatrix}0&0&0&0\\ 0&q^{-1} - q& q^2 - 1&0\\0&1 - q^{-2}& q^{-1} - q&0\\0&0&0&0\end{pmatrix},$$
which gives the normalized $R$-matrix
$$\mathcal{R}_{V_1(q^2),V_1(1)}^{norm} : V_1(1)\otimes V_1(q^2)\rightarrow V_1(q^2)\otimes V_1(1).$$
Note that this is an operator of rank $1$ which is non-invertible as opposed to $\mathcal{R}_{V,W}(z)$ 
(the reader may compare with Remark \ref{inv}).

\begin{rem} Localizing the poles of $f_{V,W}(z)\mathcal{R}_{V,W}(z)$ constitutes a difficult problem which is largely open and very important, in particular from the point of view of mathematical physics. 
\end{rem}

The Yang-Baxter equation implies that for $U$, $V$, $W$ finite-dimensional representations, the operators 
$P\circ\mathcal{R}_{.,.}(z)$ satisfy the hexagonal relation in 
$$\text{End}(U\otimes V\otimes W, W\otimes V\otimes U)[[z,w]].$$
It can be described by the following commutative diagram (the scalars are extended to the field $\mathbb{C}((z,w))$) :
\begin{equation}\label{diagYB}\xymatrix{ &V\otimes U \otimes W\ar[r]^{\text{Id}\otimes P \mathcal{R}_{U,W}(zw)}&V\otimes W\otimes U\ar[dr]^{P\mathcal{R}_{V,W}(w)\otimes\text{Id}}&   
\\ U\otimes V\otimes W \ar[ur]^{P\mathcal{R}_{U,V}(z)\otimes\text{Id}}\ar[dr]_{\text{Id}\otimes P\mathcal{R}_{V,W}(w)}& &  & W\otimes V\otimes U
\\ & U\otimes W\otimes V\ar[r]_{P\mathcal{R}_{U,W}(zw)\otimes\text{Id}}& W\otimes U\otimes V \ar[ur]_{\text{Id}\otimes P\mathcal{R}_{U,V}(z)}& }.\end{equation}
The arrows are morphisms of $\mathcal{U}_q(\hat{\mathfrak{g}})$-modules if the representation $U$ (resp. $W$) is twisted by the automorphism $\tau_z$ (resp. $\tau_{w^{-1}}$).

\begin{rem}\label{yangrem}
The quantum affine algebra $\mathcal{U}_q(\hat{\mathfrak{g}})$ has a ``rational" analogue called the ``Yangian". It is a deformation of the enveloping algebra $\mathcal{U}(\mathfrak{g}[t])$. 
It can be seen as a degeneration of a subalgebra of $\mathcal{U}_q(\hat{\mathfrak{g}})$. 
A Yangian also has a universal $R$-matrix and normalized $R$-matrices. 
\end{rem}

Other examples of $R$-matrices will also appear in the following, in the context of representations of KLR algebras (quiver Hecke algebras).

\section{Geometric construction}\label{cgeom}

Maulik-Okounkov \cite{mo} have presented a very general construction of $R$-matrices from the action of a pair of tori on a symplectic variety. It gives rise in particular to $R$-matrices which are already known, but the techniques 
which are used, ``the stable envelopes", go much further.

Nakajima varieties (or quiver varieties, see the seminar \cite{Schi}) are particularly important examples.
Indeed in a series of seminal papers Nakajima has constructed, in the equivariant $K$-theory of these varieties, certain representations of quantum affine algebras
 $\mathcal{U}_q(\hat{\mathfrak{g}})$ for $\mathfrak{g}$ simply-laced (see \cite{Nsem, Nak0} and \cite[Section 7.3]{Schi}). In the case of Yangians the corresponding representations are realized in the equivariant cohomology by Varagnolo \cite{var}. 
Moreover, the geometric study of the coproduct, by Varagnolo-Vasserot \cite{VV} and then by Nakajima \cite{Ntens}, lead to the construction of tensor products of certain finite-dimensional representations. The result of  \cite{mo} gives an
answer to a question raised naturally from these works.

\subsection{Tori and symplectic varieties}\label{toresitu}

The geometric picture is the following; consider a pair of tori
$$(\mathbb{C}^*)^r\simeq A\subset T\simeq (\mathbb{C}^*)^{r + n}$$ 
acting on an algebraic, quasi projective, smooth, symplectic variety $X$. Let $\omega\in H^0(\Omega^2X)$ be the holomorphic symplectic form of $X$. 

We assume that the line 
$$\mathbb{C} . \omega \subset H^0(\Omega^2X)$$ 
is stable under the induced action of $T$. This is a one-dimensional representation with a corresponding character. That is, a group morphism
$$h : T\rightarrow \mathbb{C}^*.$$
 Moreover we assume that 
$$A\subset \text{Ker}(h),$$ 
i.e. $\omega$ is fixed for the action of $A$.

We assume that we have a proper $T$-equivariant map
$$\pi : X\rightarrow X_0$$
for an affine $T$-variety affine $X_0$, and that $T$ is a formal $T$-variety in the sense of \cite[Section XII]{bo} 
(this is a technical condition on a spectral sequence involving the equivariant cohomology groups, see below). 
We do not necessarily assume this map $\pi$ to be a symplectic resolution of singularities.

Let $X^A$ be the set of fixed points of $X$ for the action of $A$. 

The normal fiber bundle $N(X^A)$ of $X^A$ in $X$ has a decomposition into a direct sum, indexed by characters of $A$, of bundles with a trivial action of $A$ tensored by these characters  (see \cite[Section 5.10]{CG}).

\begin{defi} The set $\Delta$ of roots of $A$ is the set of characters of $A$ occurring in the decomposition of $N(X^A)$.\end{defi}

Let $t(A)$ be the group of cocharacters of $A$ (that is of group morphisms from $\mathbb{C}^*$ to $A$). 
It is a lattice of rank $r$. We will use the space
$$\mathfrak{a}_\mathbb{R} = t(A)\otimes_\mathbb{Z} \mathbb{R} \simeq \mathbb{R}^r
.$$ 
Note that $\Delta$ can be seen as a subset of 
$$\mathfrak{a}_\mathbb{R}^* = c(A)\otimes_\mathbb{Z} \mathbb{R},$$
where $c(A)$ is the group of characters of $A$. This identification will be made in the following
without further comments.

For $\alpha\in \Delta$, we consider the orthogonal hyperplane;
$$\alpha^\perp = \{v\in \mathfrak{a}_\mathbb{R}|\alpha(v) = 0\}.$$
We get a hyperplane arrangement of $\mathfrak{a}_\mathbb{R}$ 
\begin{equation}\label{argt}\{\alpha^\perp|\alpha\in\Delta\}.\end{equation}
The complement of the union of these hyperplanes admits a partition into open subsets $\mathcal{C}_i$ called chambers which are its connected components and whose walls are the $\alpha^\perp$:
$$\mathfrak{a}_\mathbb{R}\setminus \left( \bigcup_{\alpha\in\Delta} \alpha^\perp\right) = \bigsqcup_i \mathcal{C}_i.$$
If $\sigma$ is in a chamber $\mathcal{C}_i$, is is said to be generic and  
$$X^A = X^\sigma = \{x\in X| \sigma(z).x = x,\forall z\in\mathbb{C}^*\}.$$

\begin{defi} For $\mathcal{C}$ a chamber, a point $x\in X$ is said to be $\mathcal{C}$-stable if for any $\sigma \in \mathcal{C}$ a cocharacter, the following limit exists :
$$\text{lim}_{z\rightarrow 0}(\sigma(z).x).$$
Then the limit lies in $X^A$ and does not depend on the choice of the cocharacter $\sigma\in\mathcal{C}$. We denote this limit by
$$\text{lim}_{\mathcal{C}} (x)\in X^A.$$
\end{defi}
For $Z\subset X^A$ a connected component, one defines its attracting subset
$$\text{Attr}_{\mathcal{C}}(Z) = \{x\in X\text{ $\mathcal{C}$-stable}|\text{lim}_{\mathcal{C}} (x) \in Z\}$$
which may be seen as the set of points ``attracted" by $Z$ in the spirit of the Bialynicki-Birula decomposition for projective varieties (see for example \cite[Theorem 2.4.3]{CG}).

This induces a partial ordering $\preceq_{\mathcal{C}}$ on the set of connected components of $X^A$ defined as the transitive closure of the relation : 
$$\overline{\text{Attr}_{\mathcal{C}}(Z)}\cap Z' \neq\emptyset \Rightarrow Z'\preceq_{\mathcal{C}} Z.$$
Eventually one has the full attracting subset of $Z$ which is a closed sub-variety of $X$ :
$$\text{Attr}^f_{\mathcal{C}}(Z) = \bigsqcup_{Z'\preceq Z}\text{Attr}_{\mathcal{C}}(Z').$$

For example let
$$X = T^*\mathbb{P}^n$$ 
be the cotangent bundle of $\mathbb{P}^n$ with the natural action of the torus $T = A\times \mathbb{C}^*$ where $A = (\mathbb{C}^*)^{n}$. The action of $A$ is induced from the action on $\mathbb{C}^{n+1}$ and given by characters denoted by $u_1,\ldots, u_n$ for the first $n$ coordinates. The additional factor $\mathbb{C}^*$ has a non-trivial action on the fibers of $T^*\mathbb{P}^n$ given by the character $h$. By naturally identifying $\mathbb{P}^n$ to a subvariety of  $X$  by the zero section, the fixed points in $X$ are the 
$$p_i = [0 : \ldots : 0 : 1 : 0 : \ldots : 0]$$ 
and
$$X^A = \{p_0,\ldots, p_n\}.$$ 

For example, for $n = 1$,  denote $u = u_1$. For $a\in A = \mathbb{C}^*$ and $[x,y]\in \mathbb{P}^1$, we have
$$a.[x:y] = [a^{u}x :y].$$ 
We get $2$ roots
$$\Delta=\{\alpha , -\alpha\}\subset \mathfrak{a}_{\mathbb{R}}^*\simeq \mathbb{R}$$ 
with  
$$\alpha(a) = a^{u}\text{ and }(-\alpha)(a) = a^{-u},$$ 
and $2$ chambers: 
$$\mathfrak{a}_{\mathbb{R}}  = \mathbb{R} = \mathcal{C}_+\sqcup \alpha^\perp \sqcup \mathcal{C}_-$$
where 
$$\mathcal{C}_\pm = \{x\in\mathbb{R}|\pm u.x > 0\}\text{ and }\alpha^\perp=\{0\}.$$ 
Hence 
$$\text{Attr}_{\mathcal{C}_+}(p_0) = T_{p_0}^*(\mathbb{P}^1)\text{ , }\text{Attr}_{\mathcal{C}_-}(p_0) = \mathbb{P}^1\setminus \{p_1\},$$
$$\text{Attr}_{\mathcal{C}_+}(p_1) = \mathbb{P}^1\setminus\{p_0\}\text{ , }\text{Attr}_{\mathcal{C}_-}(p_1) = T_{p_1}^*(\mathbb{P}^1).$$
In particular
$$\{p_0\}\succ_{\mathcal{C}_-}\{p_1\}\text{ et }\{p_1\}\succ_{\mathcal{C}_+}\{p_0\}.$$

\subsection{Equivariant cohomology and polarization} 

Let us briefly recall some important points about equivariant cohomology. The reader may refer to \cite{V} 
and references therein for more details. It is used as a tool to describe the topological data of a variety endowed with a
group action. It preserves the information on the stabilizers of fixed points.

For $X$ a quasi-projective smooth variety endowed with the action of a Lie group $G$, the equivariant cohomology ring is denoted by $H^\bullet_G(X)$ (we work with coefficients in $\mathbb{C}$). 
This graded $\mathbb{C}$-algebra is defined as the ordinary cohomology ring of the Borel construction (or homotopy quotient) :
$$H^{\bullet}_G(X) = H^\bullet((EG\times X)/G).$$
Here $EG$ is obtained as the total space of a universal bundle $EG\rightarrow BG$ and we consider the space of $G$-orbits on the product $EG\times X$ endowed with the diagonal group action of $G$, where $G$ acts naturally on $EG$. 
The replacement of $X$ by $EG\times X$ is justified since the action of $G$ on the product is free, and $X$ and $EG\times X$ are homotopically equivalent because $EG$ is contractible. The classifying space $BG$ 
is the homotopy quotient of a point:
$$(EG\times \{pt\})/G = BG\text{ and }H^\bullet_G(pt) = H^\bullet(BG).$$ 
This example is very important since the inclusion $\{pt\}\rightarrow X$
induces the structure of a $H^\bullet_G(pt)$-algebra on $H^\bullet_G(X)$.
In the case of a formal $G$-variety as above, one gets a free $H^\bullet_G(pt)$-module.

For $Z\subset X$ a closed subvariety stable under the action of $G$, this class $[Z]\in H^\bullet_G(X)$ is the equivariant class of its normal bundle in $X$ :
$$[Z] = e(N(Z)) \in H^\bullet_G(X).$$ 
Its degree is $codim(Z)$ (the dimensions and codimensions of varieties are given over real numbers in this section).

In the case of a torus $G =  (\mathbb{C}^*)^M$, we have 
$$H^\bullet_G(pt) = \mathbb{C}[u_1,\ldots, u_M]$$ 
where $u_1,\ldots, u_M$ are of degree $2$ and correspond to characters of $T$.

For the example $X = T^*\mathbb{P}^n$ of the previous section, we have :
$$H_T^\bullet(pt) = \mathbb{C}[u_1,\ldots, u_n,h],$$ 
$$H_T^\bullet(X) = \mathbb{C}[c,u_1,\ldots,u_n,h]/<(c - u_1)\ldots (c - u_n)c = 0>$$
where $c$ is the first Chern class of the tautological bundle $\mathcal{O}(-1)$. Hence $H_T^\bullet(X)$
is a free $H_T^\bullet(pt)$-module of rank $n + 1$ and
$$H_T^\bullet(X^A) = \bigoplus_{0\leq i\leq n} \mathbb{C}[u_1,\ldots, u_n,h] [p_i].$$

Let us go back to the general situation of Section \ref{toresitu}. 
For a chamber $\mathcal{C}$, as above, we have a decomposition
$$\Delta = \Delta_+ \sqcup  \Delta_-$$
into positive and negative roots where
$$\Delta_\pm = \{\alpha\in \Delta|\pm \alpha(\lambda) > 0, \forall \lambda \in \mathcal{C}\}.$$
Hence, for $Z\subset X^A$ a connected component, we have a decomposition :
$$N(Z) = N_+(Z) \oplus N_-(Z).$$
The symplectic form induces a duality between the two components (up to a tensor product by a trivial bundle of character $h$). In particular for $\alpha\in\Delta$ a root, we have that $-\alpha\in\Delta$ and 
$$(-1)^{\frac{\text{codim}(Z)}{2}}[Z]\in H_A^\bullet(Z)$$ 
is a perfect square in $H_A^\bullet(Z)$ :
$$\epsilon^2 = (-1)^{\frac{\text{codim}(Z)}{2}}[Z].$$ 
A polarization is a choice of a square root $\epsilon\in H_A^\bullet(Z)$. 

\subsection{Stable envelopes}

We assume a chamber $\mathcal{C}$ and a polarization have been chosen.

\begin{thm}\label{mainmo}\cite{mo} There exists a unique morphism of $H_T^\bullet(pt)$-modules 
$$\text{Stab}_{\mathcal{C}} : H_T^\bullet(X^A)\rightarrow H_T^\bullet(X)$$
such that for $\gamma\in H_{T/A}^\bullet(Z)$ where $Z\subset X^A$ is a connected component,  the image $\Gamma = \text{Stab}_{\mathcal{C}}(\gamma)$ satisfies :

(i) $\text{Supp}(\Gamma)\subset \text{Attr}^f_{\mathcal{C}}(Z)$,

(ii) $\Gamma_{\vert Z} = (\pm e(N_-(Z)))\cup\gamma$ where the sign is set so that $\pm e(N_-(Z)) = \epsilon$ in $H_A(Z)$,

(iii) for any $Z'\prec_{\mathcal{C}} Z$, the degree satisfies $\text{deg}_A (\Gamma_{\vert Z'}) < \text{codim}(Z')/2$.
\end{thm}

\noindent Note that as $A\subset T$ is commutative, the torus $T$ acts on $X^A$ with a trivial action of $A$.
Then $H_T^\bullet(X^A)\simeq H_{T/A}(X^A)\otimes_{\mathbb{C}[\mathfrak{t}/\mathfrak{a}]}\mathbb{C}[\mathfrak{t}]$ where $\mathfrak{a}\subset \mathfrak{t}$ are 
Lie algebras of $A\subset T$. The degree $\text{deg}_A$ is defined by the degree in $\mathbb{C}[\mathfrak{a}]$ 
thanks to the factorization $\mathbb{C}[\mathfrak{t}] \simeq \mathbb{C}[\mathfrak{a}]\otimes \mathbb{C}[\mathfrak{t}/\mathfrak{a}]$.

The support condition (i) means that the restriction of $\Gamma$ 
to $X\setminus\text{Attr}^f_{\mathcal{C}}(Z)$ is zero.

\noindent{\sc Proof} (sketch) --- The proof of the existence is based on a Lagrangian correspondence 
$\mathcal{L}$ defined in
$$X\times X^A\supset\mathcal{L}$$
endowed with the antidiagonal symplectic form $(\omega, - \omega_{|X_A})$. It is a stable isotropic subvariety whose dimension is half of 
the dimension of $X\times X^A$. This subvariety $\mathcal{L}$ is built by successive approximations; first one considers the closure of the
preimage of the diagonal subvariety by the natural map 
$$\text{Attr}_{\mathcal{C}}(Z)\times Z\rightarrow Z\times Z,$$
$$(x,y)\mapsto (\text{lim}_{\mathcal{C}} (x),y).$$
Then it is modified by an inductive process described in \cite[Section 3.5]{mo} so that for $Z'\succ_{\mathcal{C}} Z$, the $A$-degree of the restriction of 
$[\mathcal{L}]$ to $Z'\times Z$ is less than $\text{codim}(Z')/2$.

The reader may refer to \cite[Chapter 2]{CG} and to the seminar \cite[Sections 5.2, 7.2]{Schi} for a general discussion on 
operators defined by correspondences \footnote{The convolution product for Steinberg correspondences constitutes a very important example.}.
The construction in \cite{mo} takes place in this framework; we have the natural projections $\pi_1$ and $\pi_2$
$$\xymatrix{ & \mathcal{L}\ar[dl]_{\pi_2}\ar[dr]^{\pi_1} & \\ X^A& & X },$$
$\mathcal{L}$ is proper over $X$, and we have 
$$\xymatrix{ & H_T^\bullet(\mathcal{L})\ar[dr]^{(\pi_1)_*} & \\ H_T^\bullet(X^A)\ar[ur]_{(\pi_2)^*}& & H_T^\bullet(X)},$$
which allows us to define the map as a composition
$$\text{Stab}_{\mathcal{C}} = (\pi_1)_*\circ (\pi_2)^* : H_T^\bullet(X^A)\rightarrow H_T^\bullet(X).$$
\qed

\begin{rem}\label{isom} Note that the restriction morphism obtained from the inclusion is defined on $H_T^\bullet(X)$  :
$$H_T^\bullet(X)\hookrightarrow H_T^\bullet(X^A).$$
Under standard hypotheses, the localization theorem for equivariant cohomology implies that this morphism becomes an isomorphism
after scalar extension to the fraction field of the ring $H_T^\bullet(pt)$. The morphism $\text{Stab}_{\mathcal{C}}$ 
also becomes an isomorphism after scalar extension.
\end{rem}

\begin{rem} When $\pi$ is a symplectic resolution, the proof of Theorem \ref{mainmo} is simplified. This case includes Nakajima quiver varieties discussed below, as well as the cotangent bundle 
$T^*(G/P)$ for a parabolic subgroup $P\subset G$ of an algebraic group (with the Springer resolution).\end{rem}

Let us go back to the example $X = T^*\mathbb{P}^1$ discussed above:
$$H_T^\bullet(X) = \mathbb{C}[u,h]\oplus \mathbb{C}[u,h]c\text{ and }H_T^\bullet(X^A) = \mathbb{C}[u, h] [p_0] \oplus \mathbb{C}[u, h] [p_1] $$
(see \cite{sm} for example). The injection $i : H_T^\bullet(X)\hookrightarrow H_T^\bullet(X^A)$ corresponds to the identification:
$$[p_0] = \frac{c - u}{u}\text{ , }[p_1] = \frac{ - c}{u}.$$
 We get:
$$\text{Stab}_{\mathcal{C}_+}([p_0]) 
= u - c\text{ and }\text{Stab}_{\mathcal{C}_+}([p_1]) 
=  - c - h.$$
$$\text{Stab}_{\mathcal{C}_-}([p_0]) =u -c - h \text{ and }\text{Stab}_{\mathcal{C}_-}([p_1]) =  - c.$$
For $u\neq \pm h$, we get a basis of $H_T^\bullet(X)$ : $(\text{Stab}_{\mathcal{C}_\pm}([p_0]),\text{Stab}_{\mathcal{C}_\pm}([p_1]))$. In the basis of fixed points $([p_0], [p_1])$, we have the matrices  :
$$\text{Stab}_{\mathcal{C}_+} = \begin{pmatrix} - u & - h \\ 0 & u - h  \end{pmatrix}\text{ , }\text{Stab}_{\mathcal{C}_-} = \begin{pmatrix} - u - h&0 \\ -h& u\end{pmatrix}.$$

\begin{rem}\label{kthe} For the $K$-theoretical version of Theorem \ref{mainmo}, see \cite[Section 9]{O}. The reader may refer to \cite[Chapter 5]{CG} for an introduction 
to equivariant $K$-theory. There are differences with the cohomological case, in particular the application $\text{Stab}$ depends not only on the chamber $\mathcal{C}$, but also on an alcove 
obtained from an affine version of the arrangement (\ref{argt}). These alcoves are the connected components of
$$\mathfrak{a}_{\mathbb{R}}\setminus \left( \bigcup_{\alpha\in\Delta, n\in\mathbb{Z}} H_{\alpha,n} \right)$$
defined from affine hyperplanes
$$H_{\alpha,n}=\{v\in\mathfrak{a}_{\mathbb{R}}|\alpha(v) = n\}.$$
\end{rem}

\begin{rem} In certain cases, the applications \text{Stab} in $K$-theory have been considered by Rimanyi-Tarasov-Varchenko \cite{RTV}.
\end{rem}

\subsection{Geometric construction of $R$-matrices}

Let us fix a polarization and consider two chambers $\mathcal{C}$ and $\mathcal{C}'$. Then we have two maps $\text{Stab}_{\mathcal{C}}$ 
and $\text{Stab}_{\mathcal{C}'}$ :
$$\xymatrix{ & H_T^\bullet(X) & \\ H_T^\bullet(X^A)\ar[ur]_{\text{Stab}_{\mathcal{C}}}\ar@{-->}[rr]_{\mathcal{R}_{\mathcal{C}',\mathcal{C}}}& & H_T^\bullet(X^A)\ar[ul]^{\text{Stab}_{\mathcal{C}'}}}.$$
Up to localization at $e(N_-)$ (for $N_-$ the tautological bundle), the map $\text{Stab}_{\mathcal{C'}}$ is invertible (see Remark \ref{isom}). So after extension of scalars we can define the geometric $R$-matrix
$$R_{\mathcal{C}',\mathcal{C}} = (\text{Stab}_{\mathcal{C}'})^{-1}\circ \text{Stab}_{\mathcal{C}} \in \text{End}(H_T^\bullet(X^A)).$$
We say that one goes from the chamber $\mathcal{C}$ to the chamber $\mathcal{C}'$ by wall-crossing $\alpha^\perp$. 

For the fundamental example $X = T^*\mathbb{P}^1$ of the previous section, we recover the $R$-matrix (\ref{yang}) (The Yang $R$-matrix) :
$$R_{\mathcal{C}_-,\mathcal{C}_+} = \frac{1}{u+h}\begin{pmatrix} u & h\\ h & u\end{pmatrix}.$$

The particular case of Nakajima varieties is very important here. For $Q$ a quiver (an oriented graph) whose set of vertices is finite, 
and for $\textbf{w},\textbf{v}\in\mathbb{N}^I$, we can define the Nakajima symplectic variety\footnote{The variety depends on a stability condition which is supposed to be generic and fixed.}
$$\mathcal{M}(\textbf{w},\textbf{v}).$$

\begin{rem}\label{actiong} As we have recalled in the introduction of this section, the relation between these varieties (and their graded analogues) and quantum groups was established by Nakajima at the beginning of the theory \cite{Nsem, Nak0}. We refer to the seminar \cite{Schi} for a general introduction to these varieties and their numerous applications.
\end{rem}

Denote by $(w_i)_{i\in I}$ the coordinates of $\textbf{w}$ on $(\omega_i)_{i\in I}$ the canonical basis of $\mathbb{N}^I$ :
$$\textbf{w} = \sum_{i\in I} w_i \omega_i.$$
The variety $\mathcal{M}(\textbf{w},\textbf{v})$ is endowed with a symplectic action of the group
$$\prod_{i\in I} \text{GL}(w_i)$$
which contains a maximal torus denoted by $A = A_{\textbf{w}}$. It is contained in the larger torus
$$T = A_{\textbf{w}}\times \mathbb{C}^*,$$ 
which acts \footnote{When the quiver has loops, we can extend the action to a larger torus.} on $\mathcal{M}(\textbf{w},\textbf{v})$.
The additional factor $\mathbb{C}^*$ leads to the definition of a character $h$ as in the general context described above.

We have the following crucial factorization property. For the variety 
$$\mathcal{M}(\textbf{w}) = \bigsqcup_{\textbf{v}\in\mathbb{N}^I} \mathcal{M}(\textbf{w},\textbf{v}),$$
the space of fixed-points is itself a product of Nakajima varieties :
\begin{equation}\label{factm}(\mathcal{M}(\textbf{w}))^{A_{\textbf{w}}} = \mathcal{M}(\omega_1)^{\times w_1}\times\ldots \times \mathcal{M}(\omega_n)^{\times w_n}.\end{equation}
The reader may refer to \cite{Nt} for other recent developments related to this property. In particular
$$H_T^\bullet((\mathcal{M}(\textbf{w}))^{A_{\textbf{w}}}) = H_1^{\otimes w_1}\otimes\ldots \otimes H_n^{\otimes w_n},$$
where, for $i\in I$,
$$H_i = H_T^\bullet(\mathcal{M}(\omega_i)).$$
The following case is important. Consider
$$\textbf{w} = \omega_i + \omega_j$$ 
where $i, j\in I$. We obtain the operator
$$R_{i,j}(a)\in \text{End}(H_i\otimes H_j)\otimes \mathbb{C}(a).$$
It depends only on the difference $a = a_j - a_i$ in $\mathfrak{a}_{\mathbb{R}}$ of two cocharacters of the group $A_{\textbf{w}}\simeq (\mathbb{C}^*)^2$ associated respectively to the second factor $H_j$ and the first factor $H_i$.

\begin{prop}\cite{mo} For $i,j,k\in I$, we have in 
$$\text{End}(H_i\otimes H_j\otimes H_k,H_k\otimes H_j\otimes H_i)\otimes \mathbb{C}(a,b)$$
the following relation between operators $\mathcal{R}_{.,.}(\lambda) = P \circ R_{.,.}(\lambda)$ :
$$(\mathcal{R}_{i,j}(a)\otimes \text{Id})(\text{Id}\otimes \mathcal{R}_{i,k}(b))(\mathcal{R}_{j,k}(b-a)\otimes \text{Id}) = (\text{Id}\otimes \mathcal{R}_{j,k}(b-a))(\mathcal{R}_{i,k}(b)\otimes\text{Id}) (\text{Id}\otimes \mathcal{R}_{i,j}(a)).$$
\end{prop}

\begin{rem} The statement in this proposition means that the $R_{i,j}(a)$ satisfy the Yang-Baxter equation or that the operators $\mathcal{R}_{i,j}(a)$ make a diagram 
analogous to diagram (\ref{diagYB}) commutative. The dependence on the spectral parameter is denoted additively here.
\end{rem}

\begin{rem} In addition, the authors of \cite[Section 4.3]{mo} establish that the $R$-matrices have a remarkable factorization property in the spirit of factorizations already known for $R$-matrices 
obtained from quantum group. They can be written as a product (potentially infinite) of elementary $R$-matrices corresponding to pairs of chambers with a common wall  $\alpha^\perp$ (where $\alpha\in \Delta$ is a root).
\end{rem}

\noindent{\sc Proof} (sketch) --- Thanks to the factorization property (\ref{factm}), each factor in the relation is identified with a geometric $R$-matrix associated to a variety
$$\mathcal{M}(\omega_i + \omega_j + \omega_k)$$
for a certain chamber. The relation comes from the following general remark; consider chambers 
$$\mathcal{C} = \mathcal{C}_0,\mathcal{C}_1,\ldots ,\mathcal{C}_{2N} = \mathcal{C}$$
which contain the same $\mathcal{F}$ of codimension $2$ considered in a cyclic ordering around $\mathcal{F}$. We get
$$R_{\mathcal{C}_0,\mathcal{C}_1}\circ R_ {\mathcal{C}_1,\mathcal{C}_2}\circ \ldots \circ R_{\mathcal{C}_{2N-1},\mathcal{C}_{2N}} = \text{Id}.$$
In the situation of the statement of the proposition, with 
$$\mathcal{F} = \{x\in\mathfrak{a}_{\mathbb{R}}|a(x) = b(x) = 0\},$$ 
we have $6$ chambers $\mathcal{C}_1,\ldots,\mathcal{C}_6$ and there are two different paths from $\mathcal{C}_1$ to $\mathcal{C}_4$  :
\begin{center}
\epsfig{file=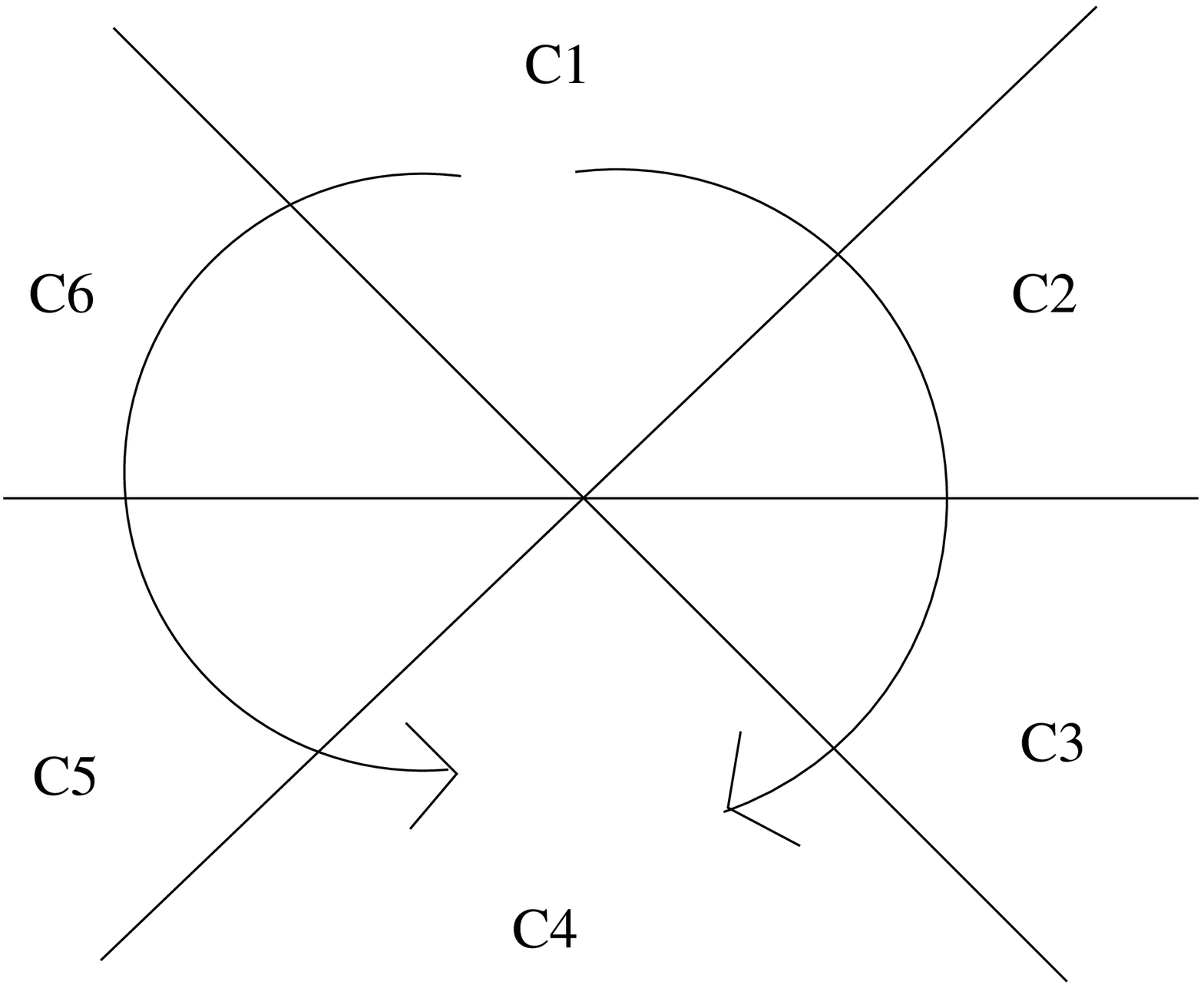,width=0.5 \linewidth}
\end{center}
We obtain the result.\qed

\subsection{Quantum groups}\label{gq}

A quantum group called the Maulik-Okounkov algebra $Y_Q$ can be constructed from $R$-matrices obtained from Nakajima varieties attached to a quiver $Q$. It is a Hopf algebra generated by all coefficients of these $R$-matrices.
More precisely, the algebra $Y_Q$ is generated by the elements of a product of algebra endomorphisms of tensor products of spaces  $H_i(u_i)$. Then each of these spaces has the structure of a representation of  $Y_Q$ and the operators  $\mathcal{R}_{i,j}(a)$ are morphisms of $Y_Q$-modules.

It is a  ``$RTT$ construction from $R$-matrices" in the spirit of \cite{FRT}: for $\mathcal{R}(z)$ the 
universal $R$-matrix of $\mathcal{U}_q(\hat{\mathfrak{g}})$ and $V$ a representation of this algebra, 
we have the images 
$$T(z) \in (\text{End}(V)\otimes \mathcal{U}_q(\hat{\mathfrak{g}}))((z)),$$
$$R(z) = \mathcal{R}_{V,V}(z) \in \text{End}(V\otimes V)((z))$$
of $\mathcal{R}(z)$. Then we have in
$$(\text{End}(V^{\otimes 2})\otimes \mathcal{U})((z,w))$$
the $RTT$ relation :
$$R_{12}(z) T_{13}(zw) T_{23}(w) = T_{23}(w) T_{13}(zw) R_{12}(z)$$
 deduced from the Yang-Baxter equation. If the  $R(z)$ are known, this $RTT$ relation can be seen as a system of relations between elements of $\mathcal{U}_q(\hat{\mathfrak{g}})$. Such a system of relations for $Y_Q$ is given by the Maulik-Okounkov construction. 

When $Q$ is a quiver of finite type, that is when the underlying diagram  $\Gamma$ is Dynkin, the algebra $Y_Q$ 
is none other than the Yangian $Y(\mathfrak{g})$ (see Remark \ref{yangrem}) associated to the Lie algebra 
$\mathfrak{g}$ of Dynkin diagram $\Gamma$ (see \cite{mc}). The action of the algebra $Y_Q$ 
on the tensor products of spaces $H_i(u_i)$ is also given by geometric constructions previously known \cite{Nsem, Nak0, var, VV}, see Remark \ref{actiong}.

\begin{rem} Let us continue the discussion of Remark \ref{kthe} : a $K$-theoretical version \cite{O} of
the construction presented in this section gives a construction of $R$-matrices and quantum affine algebras.\end{rem}

In general, it seems difficult to describe the algebra $Y_Q$. In parallel to the work of Maulik-Okounkov, Schiffmann-Vasserot \cite{SV1, SV} have associated to each quiver $Q$ an algebra called the cohomological Hall algebra. This algebra is expected to be essentially the same as the Yangian $Y_Q$. This is a topic of current study. Davison \cite{dav} 
has related these cohomological algebras to eponymous algebras introduced by Kontsevich-Soibelman.
 
The particular case of a quiver $Q$ with a single vertex and a single loop is remarkable. The corresponding algebra $Y_Q$ can be described by generators and relations \cite{mo}. It is the Yangian associated to $\hat{\mathfrak{gl}_1}$. In this case the construction in \cite{mo} is related to the construction of representations of Heisenberg algebras (see the references in \cite[Section 2.4]{Schi} and \cite{NHei} in particular).  Maulik-Okounkov \cite[Part II]{mo} and Schiffmann-Vasserot \cite{SV} have studied important applications of this Yangian $Y_Q$. The image of $Y_Q$ in
the algebra of endomorphisms of the representation $H_1^{\otimes r}$ 
can be described by using the structure of a "Verma module" of a certain affine vertex algebra 
$\mathcal{W}(\mathfrak{gl}_r)$ (in the sense of Feigin-Frenkel). This identification, more precisely the identification of a distinguished vector in the representation of $\mathcal{W}(\mathfrak{gl}_r)$ 
(the Whittaker vector) in the equivariant cohomology 
$H_1^{\otimes r}$, leads to a description of a formal power series associated to the Nakajima variety 
(the Nekrasov partition function) in terms of $\mathcal{W}(\mathfrak{gl}_r)$. This gives an answer to a question
raided by Alday-Gaiotto-Tachikawa \cite{AGT} in the framework of Conformal Field Theory (the AGT conjecture).

Other important applications of the work of Maulik-Okounkov about Baxter algebras and quantum cohomology 
will be discussed in Section \ref{cohqu}.

\section{$R$-matrices and categorifications}\label{seccat}

$R$-matrices are crucial tools used in several recent works in representation theory, in particular in the field of 
categorification of algebras.

A categorification of an algebra $\mathcal{A}$ is a monoidal category $\mathcal{M}$ endowed with 
an algebra isomorphism of its Grothendieck ring
$$\phi : \mathcal{A}\simeq K(\mathcal{M}).$$
Recall that as a group, $K(\mathcal{M})$ is the free group generated by the isomorphism classes of simple objects :
$$K(\mathcal{M}) = \bigoplus_{[V]\text{ class of a simple object in }\mathcal{M}.} \mathbb{Z} [V].$$
Then any object (not necessarily simple) in $\mathcal{M}$ has an image in $K(\mathcal{M})$ by setting  
$$[V''] = [V] + [V']$$
if $V''$ is an extension of $V$ by $V'$.
$K(\mathcal{M})$ has a ring structure by setting, for objects $V$ and $V'$ in the category $\mathcal{M}$,
$$[V\otimes V'] = [V][V'].$$

In the case of a $k$-algebra $\mathcal{A}$, $k$ a field, one can extend the scalar and replace $K(\mathcal{M})$ by $K(\mathcal{M})\otimes k$.

The notion of categorification can be refined when the algebra $\mathcal{A}$ has additional structure. 
In this case higher categorical structures are used. The reader may refer to the seminars \cite{kel, kam}. 

Two situations have been intensively studied (sometimes simultaneously) : algebras with a distinguished basis and cluster algebras.

\subsection{Categorification of algebras with basis}\label{catbas}

When the algebra $\mathcal{A}$ has a remarkable basis $\mathcal{B}$, we require the following additional condition : 

\begin{defi} A categorification $\mathcal{M}$ of an algebra $\mathcal{A}$ with a basis $\mathcal{B}$ 
is a categorification of $\mathcal{A}$ which induces a bijection between 
$\mathcal{B}$ and the basis of classes of simple objects in $K(\mathcal{M})$. 
\end{defi}

In particular this requires that the structure constants of $\mathcal{A}$ on $\mathcal{B}$ are non-negative integers.

This condition is realized in the case of the enveloping algebra $\mathcal{U}(\mathfrak{n})$ of a maximal nilpotent subalgebra  $\mathfrak{n}\subset\mathfrak{g}$ of a simple Lie algebra which is simply-laced. 
It has a remarkable basis (the canonical basis) which was constructed in particular by Lusztig as a
basis of sheaves in intersection cohomology  by interpreting the ring in terms of perverse sheaves 
(see the seminars \cite{mathieu, lit}). This construction can be seen as a categorification 
(see the seminar \cite[Section 4.1]{kam}).

A categorification in terms of representations of KLR algebras can also be used (algebras of Khovanov-Lauda-Rouquier or quiver Hecke algebras introduced in \cite{Rou, kl}). These $\mathbb{Z}$-graded algebras $R_\nu$ depend on $\mathfrak{g}$ and on a parameter $\nu\in \mathbb{N}^n$. They can be defined as sums of spaces of endomorphisms for a certain monoidal category (see the seminar \cite{kam} and references therein). In the simplest case $\mathfrak{g} = \mathfrak{sl}_2$, it is an affine nilHecke algebra. From Khovanov-Lauda \cite{kl}, Rouquier \cite{Rou} and Varagnolo-Vasserot \cite{VV2}, we get a categorification by a category of representations of KLR-algebras
$$\bigoplus_{\nu\in\mathbb{N}^n} \text{$R_\nu$-Mod}.$$ 
The product for the category is defined from a convolution product $\circ$. 
More precisely let   
$$\mathcal{U}_q(\mathfrak{n})\subset \mathcal{U}_q(\mathfrak{g})$$
be the sub-algebra corresponding to $\mathfrak{n}$ with a canonical basis compatible with the canonical basis of $\mathcal{U}(\mathfrak{n})$. The category of graded modules is a categorification of $\mathcal{U}_q(\mathfrak{n})$ 
with the basis $\mathcal{B}$ dual to the canonical basis for the natural scalar product of $\mathcal{U}_q(\mathfrak{n})$ (the canonical basis corresponds to the indecomposable projective modules). We will use the same symbol $\mathcal{B}$ for the dual canonical basis of $\mathbb{C}[N]$, the coordinate ring of a unipotent maximal subgroup $N$ of a simple complex Lie group $G$ of Lie algebra $\mathfrak{g}\supset \mathfrak{n}$. Indeed this ring $\mathbb{C}[N]$ can be seen as a specialization of $\mathcal{U}_q(\mathfrak{n})$ (different from the specialization $\mathcal{U}(\mathfrak{n})$).

We can also use a category of finite-dimensional representations of a quantum affine algebra $\mathcal{U}_q(\hat{\mathfrak{g}})$ and realize certain remarkable relations in terms of normalized $R$-matrices. As an illustration, consider the case $G = SL_3$. Then $N$ is the group of matrices $\begin{pmatrix}1 & x & z\\0 & 1 & y \\0& 0& 1\end{pmatrix}$, and the ring of polynomial functions in the coordinates $x$, $y$, $z$ is:
$$\mathbb{C}[N] = \mathbb{C}[x , y , z].$$
In this case the canonical basis is well-known and can be computed explicitly :
$$\mathcal{B} = \{x^a z^b (xy - z)^c |a,b,c\in\mathbb{N}\}\cup \{y^a z^b (xy - z)^c |a,b,c\in\mathbb{N}  \}.$$
We have a categorification of $\mathbb{C}[N]$ with its basis $\mathcal{B}$; the category $\mathcal{C}_N$ of finite-dimensional representations of $\mathcal{U}_q(\hat{\mathfrak{sl}_3})$ whose Grothendieck ring is generated by three fundamental representations 
$$[V_1(1)]\text{ , }[V_1(q^2)]\text{ , }[V_2(q)].$$ 
They correspond to $x$, $y$ and $z$, respectively.
The category $\mathcal{C}_N$ has four prime simple modules (i.e. modules which can not be factorized as a tensor product of non-trivial representations); the three fundamental representations above and another representation $W$ which corresponds to $(xy - z)$ (this module $W$ is called a Kirillov-Reshetikhin module). This can be translated as a relation in the Grothendieck ring
$$[W] = [V_1(1)][V_1(q^2)] - [V_2(q)],$$
which is obtained from an exact sequence associated to the normalized $R$-matrix $\mathcal{R}_{V_1(1),V_1(q^2)}^{norm}$
\begin{equation}\label{exacte} 0 \rightarrow W \rightarrow V_1(1)\otimes V_1(q^2) \rightarrow V_2(q)\rightarrow 0.\end{equation}
This can be seen as a particular case of results in \cite{kl, Rou, VV2} and \cite{HLClust}.

\subsection{Categorification of cluster algebras}

The reader may refer to the seminar \cite{kel} and to \cite{le} for an introduction to cluster algebras.
Fomin-Zelevinsky \cite{FZ} associated to a quiver $Q$ with $r$ vertices an algebra $\mathcal{A}(Q)$.

More precisely,  $\mathcal{A}(Q)$ is defined as a subalgebra of a fraction field
$$\mathcal{A}(Q)\subset \mathbb{Q}(X_1,\ldots, X_r).$$
A cluster is an $r+1$-tuple 
$$(z_1,\ldots, z_r, Q')$$ 
where $z_1,\ldots, z_r$ are elements of this fraction field and $Q'$ is a quiver with $r$ vertices. 
$\mathcal{A}(Q)$ is defined as the subalgebra generated by the cluster variables of all clusters obtained by iterated mutations from the initial seed 
$(X_1,\ldots, X_r, Q)$. Each cluster can be mutated in a direction $i$ where $1\leq i\leq n$ and $n\leq r$ fixed. 
In a new mutated cluster, the new variable $z_i'$ is obtained by a  Fomin-Zelevinsky mutation relation :
$$z_i' z_i =  \prod_{i\rightarrow j \text{ in }Q'} z_j  +  \prod_{j\rightarrow i\text{ in }Q'} z_j .$$
The variables $z_j$ ($j\neq i$) are not modified, but the quiver $Q'$ is replaced by the mutated quiver $\mu_i(Q')$ obtained as follows: 
for each subquiver $j\rightarrow i\rightarrow k$, an arrow $j\rightarrow k$ is added, then all arrows with source or target $i$ are reversed and
eventually all arrows of any maximal set of two by two disjoint $2$-cycles are removed.

The variables  $X_{n+1},\ldots, X_r$ which belong to all clusters are said to be frozen. A cluster monomial is a monomial in cluster variables of the same cluster.

For example, for a quiver  
$$Q : 3\longrightarrow 1\longrightarrow 2$$ 
with two frozen variables $2$ and $3$, the cluster algebra 
$$\mathcal{A}(Q)\subset \mathbb{Q}(X_1,X_2,X_3)$$ 
has $2$ clusters of cluster variables
$$(X_1,X_2,X_3)\text{ and }(X_1^{-1}(X_2 + X_3), X_2,X_3).$$ 
So there are $4$ cluster variables and one mutation relation :
\begin{equation}\label{exmuta}X_1^* X_1  = X_2 + X_3.\end{equation}
An object in a monoidal category is said to be real if its tensor square is simple.

\begin{defi} A categorification of a cluster algebra $\mathcal{A}$ is a categorification $\mathcal{M}$ 
of the algebra $\mathcal{A}$ such that the ring isomorphism
$$\phi : \mathcal{A}\simeq K(\mathcal{M})$$
sends cluster monomials of $\mathcal{A}$ to classes of real simple objects.
\end{defi}

\begin{rem}\label{forte} (i) Thanks to a factoriality property of cluster algebras \cite{glsf}, the definition implies that the cluster
variables of $\mathcal{A}$  are sent to classes of real prime simple objects.

(ii) There exists a stronger notion of categorification where it is required that  $\phi$ induces a bijection between cluster monomials
and real simple objects.

(iii) A notion of additive categorification of cluster algebras by $2$-Calabi-Yau categories has also been intensively studied, see the seminar \cite[Section 4.1]{kel} 
for an exposition of this subject.
\end{rem}

For example, the category $\mathcal{C}_N$ of representations of $\mathcal{U}_q(\hat{\mathfrak{sl}}_3)$ in the previous section is a categorification of the 
cluster algebra  $\mathcal{A}(Q)$ in the example above. Indeed there is a unique isomorphism 
$$\phi : \mathcal{A}(Q)\simeq K(\mathcal{C}_N),$$
such that
$$\phi(X_1) = [V_1(1)]\text{ , }\phi(X_2) = [V_2(q)]\text{ , }\phi(X_3) = [W]\text{ , }\phi(X_1^*) = [V_1(q^2)].$$
Then $\phi$ sends cluster variables to classes of simple objects (in fact in this case the categorification is stronger, in the sense of 
Remark \ref{forte}, (ii)). The mutation relation (\ref{exmuta}) is given by the exact sequence (\ref{exacte}) 
obtained from the normalized $R$-matrix $\mathcal{R}_{V_1(1),V_1(q^2)}^{norm}$.

Various examples of categorifications of cluster algebras have been established; first in terms of quantum affine algebras \cite{HLClust}, also in terms of perverse sheaves on Nakajima varieties \cite{N, KQ, Qin}, and of KLR algebras, which will be discussed below.

Berenstein-Zelevinsky \cite{BZ} introduced the notion of a quantum cluster algebra. These are subalgebras of quantum tori which depend on a parameter $q\in\mathbb{C}^*$. 
The specialization at $q = 1$ gives an ordinary cluster algebra. 

Geiss-Leclerc-Schr\"oer \cite{GLS} have established a fundamental example; $\mathcal{U}_q(\mathfrak{n})$ has a structure of a
quantum cluster algebra, as do certain subalgebras $A_q(\mathfrak{n}(w))$ which depend on an element  $w$ of the Weyl group of $\mathfrak{g}$ 
(they can be seen as deformations of $\mathbb{C}[N(w)]$ for a subgroup $N(w)$ of $N$). $A_q(\mathfrak{n}(w))$ can be defined as 
the space generated by the image in $A_q(\mathfrak{n})\simeq \mathcal{U}_q(\mathfrak{n})$ of products of root vectors defined by Lusztig. 
For $w = w_0$ the longest element, we recover $\mathfrak{n} = \mathfrak{n}(w_0)$.
The quantum cluster algebra structure on $A_q(\mathfrak{n}(w))$ is compatible with the 
cluster algebra structure on $\mathbb{C}[N(w)]$ (see \cite{bfz}) and in particular on $\mathbb{C}[N]$.

\subsection{$R$-matrices for KLR algebras and the Fomin-Zelevinsky conjecture}

Kang-Kashiwara-Kim-Oh established a proof of the following categorification Theorem.

\begin{thm}\label{amcate}\cite{kkko3} The cluster algebra $\mathbb{C}[N(w)]$ has a categorification.\end{thm}

\noindent{\sc Proof} (sketch) --- 
The authors use a category of representations of KLR algebras as in the results discussed in Section \ref{catbas}. 
The categorification results in \cite{kl, Rou, VV2}, as well as the theorem of Geiss-Leclerc-Schr\"oer 
on the structure of quantum cluster algebras, are crucial.
Indeed the proof relies on quantum structures. 
The corresponding notion of categorification of quantum cluster algebras was introduced and used in \cite{kkko3}.

But there is another point at the heart of the proof; the categorification of the Fomin-Zelevinsky mutation relations
in terms of exact sequences obtained from normalized $R$-matrices. For example, we have seen  (\ref{exacte}) 
an exact sequence constructed from a normalized $R$-matrix obtained from a quantum affine algebra. 
It is a mutation relation in the corresponding cluster algebra. In their proof,  Kang-Kashiwara-Kim-Oh systematically use $R$-matrices of representations of KLR algebras in order to categorify mutation relations.

For example, the exact sequence corresponding to Example (\ref{exacte}) is of the form
$$0\rightarrow L(21) \rightarrow L(1)\circ L(2) \rightarrow L(12) \rightarrow 0$$
with $L(1)$ and $L(2)$ simple $1$-dimensional modules (see \cite{KR} for complements on representations
of KLR algebras). These modules satisfy
$$L(12)\subset L(2)\circ L(1).$$

Although there is no known construction of universal $R$-matrices for KLR algebras, for $V$ and $W$ simple representations of $R_{\beta}$
and $R_{\beta'}$, respectively, there is a non-trivial morphism of $R_{\beta + \beta'}$-modules
$$\mathcal{R}^{norm}_{V,W} : V\circ W\rightarrow W\circ V.$$
Its construction relies first on the existence of intertwiners in the KLR algebra which 
lead to the definition of rational operators 
$$V\circ W\rightarrow (W\circ V)(z).$$ 
The end of the construction is analogous to the case of quantum affine algebras discussed in Section \ref{rnorm}; we have to renormalize and then to specialize.
The $R$-matrices we obtain satisfy the Yang-Baxter equation.\qed

As a consequence of Theorem \ref{amcate} one gets a proof of a general conjecture of Fomin-Zelevinsky \cite{FZ}, 
which was made precise in this case by Kimura \cite{ki} and Geiss-Leclerc-Schr\"oer \cite{GLS}. 

\begin{cor}\cite{kkko3} Cluster monomials in $\mathbb{C}[N(w)]$ are elements of the dual canonical basis.
\end{cor}

\begin{rem} (i) The case $w = w_0$ is the original motivation of the theory Fomin-Zelevinsky; to develop an
algorithmic and recursive understanding of canonical bases. These are important advances in the understanding of
these bases.

(ii) Qin \cite{Qin} has proposed another proof of this result.

(iii) As discussed in the sketch of the proof, the results above have quantum analogues when the algebras and cluster algebras are replaced
by their quantum analogues and modules over KLR algebras by their graded analogues. This is more than a generalization since the proof of Kang-Kashiwara-Kim-Oh 
uses in a crucial way the quantum structure on cluster algebras and the grading of representations of KLR algebras.
\end{rem}

\subsection{Complements : general results from intertwining operators}

Several representation theoretic results are obtained from $R$-matrices in a crucial way although these 
operators do not occur in the statements.

In particular we have the following general result (conjectured by Leclerc \cite{le1}). 
To state it we remind the reader that the head (respectively the socle) of a finite-dimensional module is its
largest semisimple quotient (resp. sub-module).

\begin{thm}\cite{kkko} Let $M$, $N$ be simple finite-dimensional modules over $\mathcal{U}_q(\hat{\mathfrak{g}})$ 
with $M$ real. Then the head and the  socle of $M\otimes N$ are real.
\end{thm}

\begin{rem} An analogous result holds for KLR algebras. It is not surprising to find analogies between quantum affine algebras and KLR algebras : results of Kang-Kashiwara-Kim \cite{kkk} give an equivalence between corresponding categories of representations. It is a generalized  Schur-Weyl duality. The construction of this duality relies on 
$R$-matrices of quantum affine algebras : their poles form a quiver which leads to the definition of
a corresponding KLR algebra.
\end{rem}

\begin{rem} Recently the proof of this statement has been adapted to the case of linear p-adic groups \cite{lm}.\end{rem}

\begin{rem} In the same spirit Kashiwara \cite{kas} proved that certain tensor products of fundamental representations are cyclic generated by a tensor product of highest weight vectors. This result has also been proved by Chari and Varagnolo-Vasserot, using other methods, for certain quantum affine algebras.\end{rem}

\section{$R$-matrices and quantum integrability}

In this last section we consider the study of quantum integrable systems which is the origin of $R$-matrices.
This is not unrelated to the work discussed in the previous two sections. In particular, the Baxter algebra (a commutative
subalgebra of a quantum group associated to a quantum system) can be interpreted in terms of quantum cohomology in
the work of Maulik-Okounkov.

\subsection{Transfer-matrices and integrability}

Let $V$ be a finite-dimensional representation of $\mathcal{U}_q(\hat{\mathfrak{g}})$ that we call 
``auxiliary space". The corresponding transfer-matrix is :
\begin{equation}\label{transfer}
\mathcal{T}_V(z) = ((\operatorname{Tr}_V \circ \rho_V) \otimes \operatorname{id})({\mathcal{R}(z)})\in \mathcal{U}_q(\hat{\mathfrak{g}})[[z]],
\end{equation}
where $\operatorname{Tr}_V$ is the trace of the space $V$ and $\rho_V$ is the morphism of the representation $V$.

We have the following fundamental ``integrability" result.

\begin{prop} Transfer-matrices commute. That is, for $V$ and $V'$ finite-dimensional representations of $\mathcal{U}_q(\hat{\mathfrak{g}})$ we have
$$\mathcal{T}_V(z)\mathcal{T}_{V'}(z') = \mathcal{T}_{V'}(z')\mathcal{T}_V(z)\text{ in }\mathcal{U}_q(\hat{\mathfrak{g}})[[z,z']].$$
\end{prop}

\noindent{\sc Proof} (sketch) --- Let us consider the image of the Yang-Baxter equation in
$$(\text{End}(V)\otimes \text{End}(V')\otimes \mathcal{U}_q(\hat{\mathfrak{g}}))[[z,z']].$$
We get a conjugation relation which gives the result after taking the trace $\text{Tr}_{(V\otimes V')}\otimes \text{Id}$.
\qed

Consider the coefficients of the transfer-matrices $\mathcal{T}_V[N]\in\mathcal{U}_q(\hat{\mathfrak{g}})$ defined by 
$$\mathcal{T}_V(z) = \sum_{N\geq 0}z^N \mathcal{T}_V[N].$$
The last Proposition means that these coefficients commute.

The properties of the trace imply that the application  $V\mapsto \mathcal{T}_V(z)$ is compatible with exact sequences. Moreover for $V$ and $V'$ in the category $\mathcal{C}$ of finite-dimensional representations of $\mathcal{U}_q(\hat{\mathfrak{g}})$ the properties of formulas (\ref{proprt}) imply
$$\mathcal{T}_{V\otimes V'}(z) = \mathcal{T}_V(z)\mathcal{T}_{V'}(z').$$
Hence we get a ring morphism 
\begin{equation}\label{take}\mathcal{T}(z) : K(\mathcal{C})\rightarrow (\mathcal{U}_q(\hat{\mathfrak{g}}))[[z]].\end{equation}
This morphism is injective and so by Frenkel-Reshetikhin \cite[Section 3.3]{Fre} we get the commutativity\footnote{The Jordan-H\"older multiplicities of simple modules are the same in $V\otimes V'$ 
and in $V'\otimes V$ but these modules are not isomorphic in general. For example for $\mathfrak{g} = sl_2$, 
$V_1(q^2)\otimes V_1(1)$ and $V_1(1)\otimes V_1(q^2)$ are not isomorphic.} 
of the Grothendieck ring $K(\mathcal{C})$.

\begin{rem}\label{trace} We also use infinite-dimensional auxiliary spaces $V$ with a natural grading 
$V = \bigoplus_{\omega = (\omega_1,\ldots,\omega_n)\in\mathbb{Z}^n} V_\omega$ by finite-dimensional weight spaces
($n$ is the rank of $\mathfrak{g}$).  Then the traces are twisted for the grading; the trace on $V$ 
is replaced by $\sum_{\omega\in\mathbb{Z}^n} (u_1^{\omega_1}\ldots u_n^{\omega_n})\text{Tr}_{V_\omega}$
for formal variables $u_i$ with $i\in I = \{1,\ldots, n\}$.\end{rem}

\begin{rem} The constructions in this subsection can be extended to the case of Yangians.\end{rem}

\subsection{Quantum systems} Transfer-matrices lead to the construction of various quantum systems; a space $W$ (the ``quantum space") on which each transfer-matrix $\mathcal{T}_V(z)$
defines an operator (the scalars are extended to $\mathbb{C}[[z]]$). 
The coefficients $\mathcal{T}_V[N]$ of the transfer-matrices act on $W$ by a large family of commutative operators.
The next step is to study the spectrum of the quantum system. That is, the eigenvalues of the transfer-matrices on $W$.

 The case $W$ a tensor product of finite-dimensional simple representations of $\mathcal{U}_q(\hat{\mathfrak{g}})$ is
a remarkable example.  Indeed for ${\mathfrak g} = \mathfrak{sl}_2$, $V = V_1$ a $2$-dimensional fundamental representation and the quantum space $W$ a tensor product of $2$ fundamental representations, 
we recover the historical case of the $XXZ$ model (see Section \ref{exfund}).
 The image of the operator $\mathcal{T}_{V_1}(z)$
in $\text{End}(W)[[z]]$ is the transfer-matrix defined by Baxter \cite{ba}, who showed\footnote{Baxter introduced the powerful method of  ``Q-operators ", which he used to solve the more involved 8 vertex model. The 6 vertex model has also been solved by other methods, in particular in the works of Lieb and Sutherland (1967).} that the eigenvalues have a remarkable form, given by the Baxter relation:
\begin{equation}\label{baxter}\lambda_j = D(z) q^{\text{deg}(Q_{j})} \frac{Q_{j}(zq^{-2})}{Q_{j}(z)} + A(z)  q^{-\text{deg}(Q_{j})} \frac{Q_{j}(zq^2)}{Q_{j}(z)},\end{equation}
where the $Q_j$ are polynomials (called Baxter polynomials) and $A(z)$, $D(z)$ are universal functions in the sense that they do not depend on the eigenvalue $\lambda_j$.

\subsection{Baxter algebras and quantum cohomology}\label{cohqu}

\begin{defi} The commutative subalgebra of $\mathcal{U}_q(\hat{\mathfrak{g}})$ generated by the coefficients $\mathcal{T}_V[N]$ of all transfer-matrices is called the Baxter algebra.
\end{defi}

The interpretation of the Baxter algebras in terms of quantum cohomology is an important application of \cite[Section 7]{mo}.

The reader may refer to  \cite{A} for a general presentation of quantum cohomology. The idea, which goes back to Vafa and Witten, is to replace the ring structure in the cohomology of a complex algebraic variety (the ordinary cup product) by a new ring structure with deformation parameters 
$q_1, q_2, \cdots, q_r$
 where $r$ is the rank of $H_2(X,\mathbb{Z})$. For $\alpha,\beta\in H^\bullet(X)$, the product in the quantum cohomology is a formal power series 
$$\alpha * \beta = \sum_{m\in\mathbb{Z}^r}(\alpha * \beta)_m q^m.$$
The coefficients $(\alpha * \beta)_m$ are defined using the theory of Gromov-Witten invariants.
We get a unitary commutative graded ring (the associativity property is remarkable) and $(\alpha * \beta)_0$ 
is the ordinary cup product.

An analogous construction can be done for equivariant cohomology. In the case of the Nakajima variety $\mathcal{M}({\bf w})$ associated to a quiver $Q$ (see Section \ref{cgeom}), for each choice of parameters $q_1,\cdots, q_r$, we get a commutative algebra of operators on $W = H_T^\bullet(\mathcal{M}({\bf w}))$ given by multiplication operators in quantum cohomology. 

Using the construction discussed in Section \ref{cgeom}, $W$ is a representation of the Maulik-Okounkov algebra $Y_Q$. The action on $W$ of the Baxter subalgebra of $Y_Q$ gives a commutative algebra of operators on $W$. But the definition of transfer-matrices can be modified by replacing the usual trace by 
$$(\text{Tr}\circ g)\otimes \text{Id}$$ 
for certain operators $g$ on the auxiliary space $V$. So Maulik-Okounkov 
consider a larger family of commutative (Baxter) algebras acting on $W$ and the actions are identified with multiplication operators in quantum cohomology.

It is a realization of the prediction made in the pioneering works \cite{NS1, NS2}.

The theory of quantum cohomology comes with a natural connection which takes into account the dependence of the product in quantum cohomology on the parameters $q_1,\cdots, q_r$. It is also related \cite{mo} to another important structure from  $R$-matrices, the quantum Knizhnik-Zamolodchikov connection in the sense of \cite{ifre}.

\subsection{$Q$-operators}

Consider the Borel subalgebra (in the sense of Drinfeld-Jimbo)
$$\mathcal{U}_q(\hat{\mathfrak{b}})\subset \mathcal{U}_q({\hat{\mathfrak g}}).$$ 
We can define the transfer-matrices associated to representations of $\mathcal{U}_q(\hat{\mathfrak{b}})$ as
$$\mathcal{R}(z)\in (\mathcal{U}_q(\hat{\mathfrak{b}})\hat{\otimes} \mathcal{U}_q({\hat{\mathfrak g}}))[[z]].$$
In order to study Baxter formulas for the $XXZ$ model, Bazhanov-Lukyanov-Zamolodchikov \cite{BLZ} constructed a family of simple representations of 
$\mathcal{U}_q(\hat{\mathfrak{b}})\subset \mathcal{U}_q(\hat{sl}_2)$ which are infinite dimensional and denoted by $L_a$ ($a\in\mathbb{C}^*$), that
we will call  ``prefundamental representations" in the rest of the text. Transfer-matrices corresponding to prefundamental representations are
called Baxter $Q$-operators (see footnote (7)).

By using the Drinfeld presentation of $\mathcal{U}_q(\hat{\mathfrak{g}})$ (see Section \ref{algcons}), we can endow $\mathcal{U}_q(\hat{\mathfrak{b}})$ itself with a triangular decomposition. That is, a linear isomorphism
$$\mathcal{U}_q(\hat{\mathfrak{b}})\simeq \mathcal{U}_q(\hat{\mathfrak{b}})^-\otimes \mathcal{U}_q(\hat{\mathfrak{b}})^0 \otimes \mathcal{U}_q(\hat{\mathfrak{b}})^+$$
with $\mathcal{U}_q(\hat{\mathfrak{b}})^-$, $\mathcal{U}_q(\hat{\mathfrak{b}})^0$, $\mathcal{U}_q(\hat{\mathfrak{b}})^+$ natural subalgebras of 
$\mathcal{U}_q(\hat{\mathfrak{b}})$. Prefundamental representations are highest weight for this triangular decomposition.

For a $2$-dimensional simple representation $V$ of $\mathcal{U}_q(\hat{sl}_2)$, we have 
\begin{equation}\label{exexa}0 \rightarrow  [\omega]\otimes L_{q^{-2}} \rightarrow V\otimes L_1 \rightarrow [-\omega] \otimes L_{q^2} \rightarrow 0,\end{equation}
a non-split exact sequence, with $[\pm \omega]$ invertible representations of dimension $1$. Taking transfer-matrices, we recover Formula (\ref{baxter}).

\begin{rem} For the $8$-vertex model, a construction of $Q$-operators in terms of representations of elliptic quantum groups was discovered recently by Felder-Zhang \cite{FZ2}.
\end{rem}

\subsection{Frenkel-Reshetikhin Approach}

Frenkel-Reshetikhin \cite{Fre} proposed\footnote{An analogous conjecture had been formulated in particular cases by Reshetikhin, Bazhanov-Reshetikhin and Kuniba-Suzuki.} a new approach in order to generalize Baxter formulas. In this framework they introduce the $q$-character $\chi_q(V)$ of a finite-dimensional representation $V$ of $\mathcal{U}_q(\hat{\mathfrak{g}})$. It is a Laurent polynomial in the indeterminates $Y_{i,a}$ ($1\leq i\leq n$, $a\in\mathbb{C}^*$) 
$$\chi_q(V) \in \mathbb{Z}[Y_{i,a}^{\pm 1}]_{1\leq i\leq n, a\in\mathbb{C}^*}.$$
 The definition of $\chi_q(V)$ relies on the decomposition of $V$ into Jordan subspaces (for a commutative family of operators on $W$ obtained from the Drinfeld realization). The coefficients of $\chi_q(V)$ are positive and the sum is $\text{dim}(V)$. 

For example, for $\mathfrak{g} = \mathfrak{sl}_2$ and $V = V_1$ a fundamental representation of dimension $2$,  
\begin{equation}\label{qcar}
\chi_q(V) = Y_{1,q^{-1}} + Y_{1,q}^{-1}.
\end{equation}
For $1\leq i\leq n$, we set $q_i = q^{r_i}$ where $r_i$ is the length of the corresponding simple root.

\begin{conj}\label{tconj}\cite[Section 6.1]{Fre} For $V$ a finite-dimensional representation, the eigenvalues 
$\lambda_j$ of $\mathcal{T}_V(z)$ on a tensor product $W$ of finite-dimensional simple representations are obtained in the following way : we replace in $\chi_q(V)$ each variable $Y_{i,a}$ by
$$F_{i}(az) q^{\text{deg}(Q_{i,j})} \frac{Q_{i,j}(zaq_i^{-1})}{Q_{i,j}(zaq_i)},$$
where the $Q_{i,j}(z)$ are polynomials and $F_{i}(z)$ does not depend on the eigenvalue $\lambda_j$.
\end{conj}

The $Q_{i,j}(z)$ are generalized Baxter polynomials. It is a generalization of the results of Baxter: by using this process we recover formula (\ref{baxter}) from formula (\ref{qcar}). In the general case, the generalized Baxter formula for the $\lambda_j$ has $\text{dim}(V)$ terms.

The statement of Conjecture \ref{tconj} with the deformed trace (as in Remark \ref{trace}) is established in \cite{FH}.

\noindent{\sc Proof} (sketch) --- The idea is to interpret the $Q_i$ as eigenvalues of new
transfer-matrices, constructed from prefundamental representations $L_{i,a}$ (where $1\leq i\leq n$ and $a\in\mathbb{C}^*$) of $\mathcal{U}_q(\hat{\mathfrak{b}})$ \cite{HJ} generalizing the representations of \cite{BLZ} discussed above.
It is proved that, up to a universal scalar factor $F_i(z)$, the transfer-matrix associated to $\mathcal{T}_i(z)$ 
acts on the quantum space $W$ by a polynomial operator. Then we have to prove that the eigenvalues of 
$\mathcal{T}_V(z)$ can be written in terms of the eigenvalues of the $\mathcal{T}_i(z)$ using the $q$-character of $V$. 
To do this, one can use a category $\mathcal{O}$ of representations of $\mathcal{U}_q(\hat{\mathfrak{b}})$ containing the finite-dimensional representations as well as the prefundamental representations. Consider its Grothendieck ring\footnote{As the representations in the category $\mathcal{O}$ are not of finite length, the multiplicity of a simple module in a module is defined as in \cite{ka}, see \cite[Section 3.2]{HL3}.} $K(\mathcal{O})$.
By substituting in $\chi_q(V)$ each $Y_{i,a}$ by\footnote{Here $[\omega_i]$ is the class of a $1$-dimensional representation.} $[\omega_i]\frac{[L_{i,aq^{-1}}]}{[L_{i,aq}]}$ and $\chi_q(V)$ by $[V]$, one gets a relation 
in the fraction field of $K(\mathcal{O})$ (this can be compared to the example in formula (\ref{qcar}) 
with the exact sequence (\ref{exexa})).
Now we can extend the transfer-matrix map (\ref{take}) to $K(\mathcal{O})$. Then the relations in $K(\mathcal{O})$ imply the desired relations between transfer-matrices and eigenvalues.
\qed

\subsection{Bethe Ansatz Equations}
It is quite easy to write an explicit formula for the scalar functions $F_i(z)$. It is much more complicated to get 
information on the polynomial part $(F_i(z))^{-1}\mathcal{T}_i(z)$. Its eigenvalues are the Baxter polynomials $(Q_i(z))_{i\in I}$. 
The Bethe Ansatz equations between the roots of these polynomials have been conjectured (see \cite{R3, Fre})  : 

\begin{equation}\label{bethe}v_i^{-1} \prod_{j \in I}
\frac{{Q}_j(wq^{B_{ij}})}{{ Q}_j(wq^{-B_{ij}})} = -1,\end{equation}
for $i\in I$ and $w$ a root of $Q_i$. Here the $B_{i,j} = (\alpha_i,\alpha_j)$ are the coefficients of the symmetrized Cartan 
matrix and $v_i = \prod_{j\in I}u_i^{C_{j,i}}$.

\begin{rem}\label{rmref} The equations (\ref{bethe}) have been established recently under generecity conditions 
\cite{FH2} thanks to relations in $K(\mathcal{O})$ called $Q\tilde{Q}$-systems (discovered \cite{mrv} in the Langlands dual context, see above), then without generecity condition \cite{FJMM}, 
thanks in particular to a two term version of the Baxter formulas (the $QQ^*$-systems in \cite{HL3} which occur naturally as mutation relations in the spirit of the results in Section \ref{seccat}).\end{rem}

\begin{rem} Very recently \cite{PSZ}, the Baxter operators for $\mathcal{U}_q(\hat{\mathfrak{sl}_2})$ have been realized in the framework of the Maulik-Okounkov theory (Section \ref{cgeom}), 
in the form of multiplication in quantum cohomology (resp. $K$-theory) (see Section \ref{cohqu}). The Bethe Ansatz equations (\ref{bethe}) occur naturally in this context.
\end{rem}

\subsection{Other quantum systems and Langlands duality} The ODE/IM correspondence (for ``Ordinary Differential Equation" and ``Integrable model")
discovered in \cite{DT, BLZ2} gives a surprising link between, on one hand, functions associated to Schr\"odinger differential operators and, on
the other hand, the eigenvalues of a quantum system called ``quantum KdV". The latter have been constructed \cite{BLZ2} using transfer-matrices $\mathcal{T}_V(z)$ 
acting no longer on a representation of $\mathcal{U}_q(\hat{\mathfrak{g}})$ but on Fock spaces of a quantum Heisenberg algebra. 
The differential operators are of the form 
$-\partial^2 + x^{2M} + \frac{\ell(\ell - 1)}{x^2}$ with  $M > 0$ integer and $\ell\in\mathbb{C}$. The corresponding functions are the spectral determinants defined as 
coefficients between bases of natural solutions. 
It is remarkable that they satisfy the Baxter relations (\ref{baxter}), although they come from a context which is {\it a priori} very different.

Feigin-Frenkel \cite{FF} have described this correspondence as an instance of the Langlands duality. The Schr\"odinger differential operators are generalized by ``affine opers" which are differential systems associated to the Langlands dual Lie algebra of $\mathfrak{g}$. The conjectures stated in this paper \cite{FF} are largely open but it is a fruitful course of 
inspiration in the study of representations of $\mathcal{U}_q(\hat{\mathfrak{g}})$ (see for example the works in \cite{mrv} mentioned in Remark \ref{rmref}).

\bigskip
\bigskip

{\bf Acknowledgments}

\medskip

I thank warmly Caroline Gruson, Bernard Leclerc, Hiraku Nakajima, Andrei Okounkov, Marc Rosso, Olivier Schiffmann, Eric Vasserot and Mich\`ele Vergne for reading the text and for their precious advise, as well as Masaki Kashiwara, Myungho Kim and Gus Schrader, whose recent mini-courses were very useful. 
Thanks also to the participants of the conference " Geometry and Representation Theory" at the Institut Erwin Schr\"odinger in Vienna where the seminar 
talk has been tested.

Special thanks to Ruari Walker for his help with the English translation.

While this text was being written, Maulik-Okounkov made public a new version of their book \cite{mo} to which we refer.

\end{document}